\documentclass[12pt]{amsart}
\usepackage{amssymb}
\usepackage{t1enc}
\usepackage[latin2]{inputenc}
\usepackage{verbatim}
\usepackage{amsmath,amsfonts,amssymb,amsthm}
\usepackage[mathcal]{eucal}
\usepackage{enumerate}
\usepackage[centertags]{amsmath}
\usepackage{graphics}

\newtheorem{theorem}{Theorem}
\newtheorem{lemma}{Lemma}

\numberwithin{equation}{subsection}

\begin{document}
\author{George Tephnadze}
\title[Reisz logarithmic means ]{On the maximal operators of Riesz
logarithmic means of Vilenkin-Fourier series}
\address{G. Tephnadze, Department of Mathematics, Faculty of Exact and
Natural Sciences, Tbilisi State University, Chavchavadze str. 1, Tbilisi
0128, Georgia and Department of Engineering Sciences and Mathematics, Lule%
\aa {} University of Technology, SE-971 87, Lule\aa {}, Sweden}
\email{giorgitephnadze@gmail.com}
\thanks{The research was supported by Shota Rustaveli National Science
Foundation grant no.13/06 (Geometry of function spaces, interpolation and
embedding theorems}
\date{}
\maketitle

\begin{abstract}
The main aim of this paper is to investigate $\left( H_{p},L_{p}\right) $
and $\left( H_{p},L_{p,\infty }\right) $ type inequalities for maximal
operators of Riesz logarithmic means of one-dimensional Vilenkin-Fourier
series.
\end{abstract}

\date{}



\section{ Introduction}

Weak (1,1)-type inequality for the maximal operator of Fejér means $\sigma
^{\ast }$ for Walsh-Fourier series was proved by Schipp \cite{Sc} and for
Vilenkin system by Pál, Simon \cite{PS}. Fujji \cite{Fu} and Simon \cite{Si2}
verified that the $\sigma ^{\ast }$ is bounded from $H_{1}$ to $L_{1}$.
Weisz \cite{We2} generalized this result and proved the boundedness of $%
\sigma ^{\ast }$ from the martingale Hardy space $H_{p}$ to the space $%
L_{p}, $ for $p>1/2$. Simon \cite{Si1} gave a counterexample, which shows
that boundedness does not hold for $0<p<1/2.$ The counterexample for $p=1/2$
due to Goginava \cite{BGG2}, (see also \cite{Go} and \cite{tep1}).

Weisz \cite{We3} proved that following is true:

\textbf{Theorem W. }The maximal operator of Fejér means $\sigma ^{\ast }$ is
bounded from the Hardy space $H_{1/2}$ to the space $L_{1/2,\infty }.$

In \cite{tep2} and \cite{tep3} it were proved that the maximal operator $%
\widetilde{\sigma }$ $_{p}^{\ast \,},$ defined by
\begin{equation*}
\widetilde{\sigma }_{p}^{\ast }:=\sup_{n\in \mathbb{N}}\frac{\left\vert
\sigma _{n}\right\vert }{\left( n+1\right) ^{1/p-2}\log ^{2\left[ 1/2+p%
\right] }\left( n+1\right) },
\end{equation*}%
where $0<p\leq 1/2$ and $\left[ 1/2+p\right] $ denotes integer part of $%
1/2+p,$ is bounded from the Hardy space $H_{p}$ to the space $L_{p}.$

Moreover, for any nondecreasing function $\varphi :\mathbb{N}_{+}\rightarrow
\lbrack 1,$ $\infty )$ satisfying the condition
\begin{equation}
\overline{\lim_{n\rightarrow \infty }}\frac{\left( n+1\right) ^{1/p-2}\log
^{2\left[ 1/2+p\right] }\left( n+1\right) }{\varphi \left( n\right) }%
=+\infty ,  \label{cond1}
\end{equation}%
there exists a martingale $f\in H_{p},$ such that
\begin{equation*}
\underset{n\in \mathbb{N}}{\sup }\left\Vert \frac{\sigma _{n}f}{\varphi
\left( n\right) }\right\Vert _{p}=\infty .
\end{equation*}

For Walsh-Paley system analogical theorem is proved in \cite{GoSzeged} and
for Walsh-Kaczmarz system in \cite{GNCz} and \cite{tep5}.

Riesz` s logarithmic means with respect to the Walsh system was studied by
Simon \cite{Si1}, Goginava \cite{gogin}, Gát, Nagy \cite{GN} and for
Vilenkin systems by Gát \cite{Ga1}, Blahota, Gát \cite{bg}, Tephnadze \cite%
{tep4}. In this paper it was proved that maximal operator of Riesz
logarithmic means of Vilenkin-Fourier series is bounded from the martingale
Hardy space $H_{p}$ to the space $L_{p}$ when $p>1/2$ and is not bounded
from the martingale Hardy space $H_{p}$ to the space $L_{p}$ when $0<p\leq
1/2.$

The main aim of this paper is to investigate $\left( H_{p},L_{p}\right) $
and $\left( H_{p},L_{p,\infty }\right) $ type inequalities for weighted
maximal operators of Riesz logarithmic means of one-dimensional
Vilenkin-Fourier series.

\section{Definitions and Notations}

Let $\mathbf{P}_{+}$ denote the set of the positive integers , $\mathbf{P}:=%
\mathbf{P}_{+}\cup \{0\}.$

Let $m:=(m_{0,}m_{1},...)$ denote a sequence of the positive integers not
less than 2.

Denote by
\begin{equation*}
Z_{m_{k}}:=\{0,1,...m_{k}-1\}
\end{equation*}
the additive group of integers modulo $m_{k}.$

Define the group $G_{m}$ as the complete direct product of the group $%
Z_{m_{j}}$ with the product of the discrete topologies of $Z_{m_{j}}$ $^{,}$%
s.

The direct product $\mu $ of the measures
\begin{equation*}
\mu _{k}\left( \{j\}\right) :=1/m_{k}\text{ \qquad }(j\in Z_{m_{k}})
\end{equation*}%
is a Haar measure on $G_{m_{\text{ }}}$with $\mu \left( G_{m}\right) =1.$

If $\sup_{n}m_{n}<\infty $, then we call $G_{m}$ a bounded Vilenkin group.
If the generating sequence $m$ is not bounded then $G_{m}$ is said to be an
unbounded Vilenkin group. \textbf{In this paper we discuss bounded Vilenkin
groups only.}

The elements of $G_{m}$ are represented by sequences
\begin{equation*}
x:=(x_{0},x_{1},...,x_{j},...)\qquad \left( \text{ }x_{k}\in
Z_{m_{k}}\right) .
\end{equation*}

It is easy to give a base for the neighborhood of $G_{m}$
\begin{equation*}
I_{0}\left( x\right) :=G_{m},
\end{equation*}%
\begin{equation*}
I_{n}(x):=\{y\in G_{m}\mid y_{0}=x_{0},...y_{n-1}=x_{n-1}\},\text{ \ }(x\in
G_{m},\text{ }n\in \mathbf{P}).
\end{equation*}

Denote $I_{n}:=I_{n}\left( 0\right) $ for $n\in \mathbf{P}$ and $\overline{%
I_{n}}:=G_{m}$ $\backslash $ $I_{n}$.

Let

\begin{equation*}
e_{n}:=\left( 0,0,...,x_{n}=1,0,...\right) \in G_{m}\qquad \left( n\in
\mathbf{P}\right) .
\end{equation*}

It is evident
\begin{equation}
\overline{I_{M}}=\left( \overset{M-2}{\underset{k=0}{\bigcup }}\overset{%
m_{k}-1}{\underset{x_{k}=1}{\bigcup }}\overset{M-1}{\underset{l=k+1}{\bigcup
}}\overset{m_{l}-1}{\underset{x_{l}=1}{\bigcup }}I_{l+1}\left(
x_{k}e_{k}+x_{l}e_{l}\right) \right) \bigcup \left( \underset{k=1}{%
\bigcup\limits^{M-1}}\overset{m_{k}-1}{\underset{x_{k}=1}{\bigcup }}%
I_{M}\left( x_{k}e_{k}\right) \right) .  \label{3}
\end{equation}

If we define the so-called generalized number system based on $m$ in the
following way :
\begin{equation*}
M_{0}:=1,\text{ \qquad }M_{k+1}:=m_{k}M_{k\text{ }}\ \qquad (k\in \mathbf{P}%
).
\end{equation*}
then every $n\in \mathbf{P}$ can be uniquely expressed as $n=\overset{\infty
}{\underset{k=0}{\sum }}n_{j}M_{j}$ where $n_{j}\in Z_{m_{j}}$ $~(j\in
\mathbf{P})$ and only a finite number of $n_{j}`$s differ from zero. Let $%
\left| n\right| :=\max $ $\{j\in \mathbf{P};$ $n_{j}\neq 0\}.$

Denote by $L_{1}\left( G_{m}\right) $ the usual (one dimensional) Lebesgue
space.

Next, we introduce on $G_{m}$ an orthonormal system which is called the
Vilenkin system.

At first, define the complex valued function $r_{k}\left( x\right)
:G_{m}\rightarrow
\mathbb{C}
,$ the generalized Rademacher functions as
\begin{equation*}
r_{k}\left( x\right) :=\exp \left( 2\pi ix_{k}/m_{k}\right) \text{ \qquad }%
\left( i^{2}=-1,\text{ }x\in G_{m},\text{ }k\in \mathbf{P}\right) .
\end{equation*}

Now, define the Vilenkin system $\psi :=(\psi _{n}:n\in \mathbf{P})$ on $%
G_{m}$ as:
\begin{equation*}
\psi _{n}(x):=\overset{\infty }{\underset{k=0}{\Pi }}r_{k}^{n_{k}}\left(
x\right) \text{ \qquad }\left( n\in \mathbf{P}\right) .
\end{equation*}

Specifically, we call this system the Walsh-Paley one if $m\equiv 2$.

The Vilenkin system is orthonormal and complete in $L_{2}\left( G_{m}\right)
$ \cite{AVD}.

Now, we introduce analogues of the usual definitions in Fourier-analysis.

If $f\in L_{1}\left( G_{m}\right) $ we can establish the Fourier
coefficients, the partial sums of the Fourier series, the Fejér means, the
Dirichlet and Fejér kernels with respect to the Vilenkin system $\psi $ in
the usual manner:%
\begin{eqnarray*}
\widehat{f}\left( k\right) &:&=\int_{G_{m}}f\overline{\psi }_{k}d\mu ,\text{
\ \ \ \ \ }\left( \text{ }k\in \mathbf{P}\text{ }\right) , \\
S_{n}f &:&=\sum_{k=0}^{n-1}\widehat{f}\left( k\right) \psi _{k}\ ,\text{ \ \
\ }\left( \text{ }n\in \mathbf{P}_{+},\text{ }S_{0}f:=0\text{ }\right) , \\
\sigma _{n}f &:&=\frac{1}{n}\sum_{k=0}^{n-1}S_{k}f\text{ },\text{ \ \ \ \ \
\ }\left( \text{ }n\in \mathbf{P}_{+}\text{ }\right) , \\
D_{n} &:&=\sum_{k=0}^{n-1}\psi _{k\text{ }},\text{ \qquad\ \ \ \ \ }\left(
\text{ }n\in \mathbf{P}_{+}\text{ }\right) , \\
K_{n} &:&=\frac{1}{n}\overset{n-1}{\underset{k=0}{\sum }}D_{k},\text{ \ \ \
\ \ \ }\left( \text{ }n\in \mathbf{P}_{+}\text{ }\right) .
\end{eqnarray*}%
$\qquad $ $\qquad $

Recall that
\begin{equation}
\quad \hspace*{0in}D_{M_{n}}\left( x\right) =\left\{
\begin{array}{l}
\text{ }M_{n},\text{\thinspace \thinspace \thinspace \thinspace if\thinspace
\thinspace }x\in I_{n}, \\
\text{ }0,\text{\thinspace \thinspace \thinspace \thinspace \thinspace if
\thinspace \thinspace }x\notin I_{n}.%
\end{array}
\right.  \label{4}
\end{equation}

It is well-known that \vspace{0pt}

\begin{equation}
\sup_{n}\int_{G_{m}}\left\vert K_{n}\right\vert d\mu \leq c<\infty .
\label{6}
\end{equation}

The norm (or quasi-norm) of the space $L_{p}(G_{m})$ is defined by \qquad
\qquad \thinspace \
\begin{equation*}
\left\| f\right\| _{p}:=\left( \int_{G_{m}}\left| f\right| ^{p}d\mu \right)
^{1/p}\qquad \left( 0<p<\infty \right) .
\end{equation*}

The space $L_{p,\infty }\left( G\right) $ consists of all measurable
functions $f$ for which

\begin{equation*}
\left\Vert f\right\Vert _{L_{p,\infty }(G)}:=\underset{\lambda >0}{\sup }%
\lambda ^{p}\mu \left( f>\lambda \right) <+\infty .
\end{equation*}

The $\sigma $-algebra generated by the intervals $\left\{ I_{n}\left(
x\right) :x\in G_{m}\right\} $ will be denoted by $\digamma _{n}$ $\left(
n\in \mathbf{P}\right) .$ Denote by $f=\left( f^{\left( n\right) },n\in
\mathbf{P}\right) $ a martingale with respect to $\digamma _{n}$ $\left(
n\in \mathbf{P}\right) .$ (for details see e.g. \cite{We1}). The maximal
function of a martingale $f$ is defend by \qquad
\begin{equation*}
f^{\ast }=\sup_{n\in \mathbf{P}}\left\vert f^{\left( n\right) }\right\vert ,
\end{equation*}%
respectively.

In case $f\in L_{1},$ the maximal functions are also given by
\begin{equation*}
f^{\ast }\left( x\right) =\sup_{n\in \mathbf{P}}\frac{1}{\left\vert
I_{n}\left( x\right) \right\vert }\left\vert \int_{I_{n}\left( x\right)
}f\left( u\right) \mu \left( u\right) \right\vert .
\end{equation*}

For $0<p<\infty $ the Hardy martingale spaces $H_{p}$ $\left( G_{m}\right) $
consist of all martingales for which
\begin{equation*}
\left\Vert f\right\Vert _{H_{p}}:=\left\Vert f^{\ast }\right\Vert
_{p}<\infty .
\end{equation*}

If $f\in L_{1},$ then it is easy to show that the sequence $\left(
S_{M_{n}}\left( f\right) :n\in \mathbf{P}\right) $ is a martingale. If $%
f=\left( f^{\left( n\right) },n\in \mathbf{P}\right) $ is martingale then
the Vilenkin-Fourier coefficients must be defined in a slightly different
manner: $\qquad \qquad $
\begin{equation*}
\widehat{f}\left( i\right) :=\lim_{k\rightarrow \infty
}\int_{G_{m}}f^{\left( k\right) }\left( x\right) \overline{\psi }_{i}\left(
x\right) d\mu \left( x\right) .
\end{equation*}%
\qquad \qquad \qquad \qquad

The Vilenkin-Fourier coefficients of $f\in L_{1}\left( G_{m}\right) $ are
the same as those of the martingale $\left( S_{M_{n}}\left( f\right) :n\in
\mathbf{P}\right) $ obtained from $f$ .

In the literature, there is the notion of Riesz` s logarithmic means of the
Fourier series. The $n$-th Riesz`s logarithmic means of the Fourier series
of an integrable function $f$ is defined by

\begin{equation*}
R_{n}f:=\frac{1}{l_{n}}\overset{n}{\underset{k=1}{\sum }}\frac{S_{k}f}{k},
\end{equation*}%
where $l_{n}:=\sum_{k=1}^{n}\frac{1}{k}.$

The kernels of Riesz`s logarithmic means is established by

\begin{equation*}
L_{n}:=\frac{1}{l_{n}}\overset{n}{\underset{k=1}{\sum }}\frac{D_{k}\left(
x\right) }{k}.
\end{equation*}

For the martingale $f$ we consider the following maximal operators

\begin{eqnarray*}
\sigma ^{\ast }f &:&=\sup_{n\in \mathbf{P}}\left\vert \sigma
_{n}f\right\vert ,\text{ \ \ \ \ \ \ \ \ \ \ \ \ \ }R^{\ast }f:=\underset{%
n\in \mathbf{P}}{\sup }\left\vert R_{n}f\right\vert , \\
\overset{\sim }{R}^{\ast }f &:&=\underset{n\in \mathbf{P}}{\sup }\frac{%
\left\vert R_{n}f\right\vert }{\log \left( n+1\right) },\text{ \ }\overset{%
\sim }{R}_{p}^{\ast }f:=\underset{n\in \mathbf{P}}{\sup }\frac{\log \left(
n+1\right) \left\vert R_{n}f\right\vert }{\left( n+1\right) ^{1/p-2}}.
\end{eqnarray*}%
A bounded measurable function $a$ is $p$-atom, if there exist a dyadic
interval $I$, such that \qquad
\begin{equation*}
\int_{I}ad\mu =0,\text{ \ }\left\Vert a\right\Vert _{\infty }\leq \mu \left(
I\right) ^{-1/p},\text{\ \ supp}\left( a\right) \subset I.
\end{equation*}

\section{Formulation of Main Results}

\begin{theorem}
The maximal operator of Riesz logarithmic means $R^{\ast }$ is bounded from
the Hardy space $H_{1/2\text{ }}$ to the space $L_{1/2,\infty }.$
\end{theorem}

Earlier, It was proved that the maximal operator $R^{\ast }$ is not bounded
form the the Hardy space $H_{1/2\text{ }}$ to the space $L_{1/2}.$ So, it is
interesting to discuss that what type weight we have to apply to get back
the boundedness of the maximal operator. We found the answer in the next
theorem.

\begin{theorem}
a)\bigskip\ The maximal operator$\ \overset{\sim }{R}^{\ast }$\textit{is
bounded from the Hardy space }$H_{1/2}$\textit{\ to the space }$L_{1/2}.$
\end{theorem}

\textit{b) Let} $\varphi :\mathbf{P}_{+}\rightarrow \lbrack 1,$\textit{\ }$%
\infty )$\textit{\ be a nondecreasing function satisfying the condition}

\begin{equation}
\overline{\lim_{n\rightarrow \infty }}\frac{\log \left( n+1\right) }{\varphi
\left( n\right) }=+\infty .  \label{55}
\end{equation}%
\textit{Then the maximal operator}

\begin{equation*}
\sup_{n\in \mathbf{P}}\frac{\left\vert R_{n}f\right\vert }{\varphi \left(
n\right) }
\end{equation*}%
\textit{is not bounded from the Hardy space }$H_{1/2}$\textit{\ to the space
}$L_{1/2}.$

\begin{theorem}
a) \bigskip Let $0<p<1/2.$ Then the maximal operator $\overset{\sim }{R}%
_{p}^{\ast }$ \textit{is bounded from the Hardy space }$H_{p}$\textit{\ to
the space }$L_{p}.$
\end{theorem}

\textit{b) Let }$0<p<1/2$\textit{\ and }$\varphi :\mathbf{P}_{+}\rightarrow
\lbrack 1,$\textit{\ }$\infty )$\textit{\ be a nondecreasing function
satisfying the condition}

\begin{equation}
\frac{\left( n+1\right) ^{1/p-2}}{\log \left( n+1\right) \varphi \left(
n\right) }=\infty .  \label{66}
\end{equation}%
\textit{Then the maximal operator}

\begin{equation*}
\sup_{n\in \mathbf{P}}\frac{\left\vert R_{n}f\right\vert }{\varphi \left(
n\right) }
\end{equation*}%
\textit{is not bounded from the Hardy space }$H_{p}$\textit{\ to the space }$%
L_{p,\infty }.$

\section{AUXILIARY PROPOSITIONS}

\begin{lemma}
\cite{We4} (Weisz) A martingale $f=\left( f^{\left( n\right) },n\in \mathbf{P%
}\right) $ is in $H_{p}\left( 0<p\leq 1\right) $ if and only if there exist
a sequence $\left( a_{k},k\in \mathbf{P}\right) $ of p-atoms and a sequence $%
\left( \mu _{k},k\in \mathbf{P}\right) $ of a real numbers such that for
every $n\in P$
\end{lemma}

\begin{equation}
\qquad \sum_{k=0}^{\infty }\mu _{k}S_{M_{n}}a_{k}=f^{\left( n\right) },
\label{6.1}
\end{equation}

\begin{equation*}
\qquad \sum_{k=0}^{\infty }\left\vert \mu _{k}\right\vert ^{p}<\infty .
\end{equation*}%
\textit{Moreover, }$\left\Vert f\right\Vert _{H_{p}}\backsim \inf \left(
\sum_{k=0}^{\infty }\left\vert \mu _{k}\right\vert ^{p}\right) ^{1/p}$%
\textit{, where the infimum is taken over all decomposition of }$f$\textit{\
of the form (\ref{6.1}).}

\begin{lemma}
\cite{Gat} (Gát) Let $A>t,$ $t,A\in \mathbf{P},$ $x\in I_{t}\backslash $ $%
I_{t+1}.$ Then
\end{lemma}

\begin{equation*}
K_{2^{A}}\left( x\right) =\left\{
\begin{array}{c}
\text{ }2^{t-1},\text{\ \ \ \ \ \ \ \ \ \ \ \ \ \ \ \ \ \ if \ \ }x\in
I_{A}\left( e_{t}\right) , \\
\left( 2^{A}+1\right) /2,\text{ \ if \ \ }x\in I_{A},\text{\ } \\
0,\text{\qquad\ \ \ \ \ \ \ \ \ \ otherwise.\ \ }%
\end{array}%
\right.
\end{equation*}

Analogously of Lemma 4 in \cite{tep3} if we apply Lemma 2 we can prove that
following is true:

\begin{lemma}
Let $x\in I_{N}\left( x_{k}e_{k}+x_{l}e_{l}\right) ,$ $\ 1\leq x_{k}\leq
m_{k}-1,$ $1\leq x_{l}\leq m_{l}-1,$ $k=0,...,N-2,$ $\ l=k+1,...,N-1.$ Then
\end{lemma}

\begin{equation*}
\int_{I_{N}}\left| K_{n}\left( x-t\right) \right| d\mu \left( t\right) \leq
\frac{cM_{l}M_{k}}{nM_{N}},\qquad \text{when }n\geq M_{N}.
\end{equation*}

Let $x\in I_{N}\left( x_{k}e_{k}\right) ,$ $1\leq x_{k}\leq m_{k}-1,$ $%
k=0,...,N-1.$ Then

\begin{equation*}
\int_{I_{N}}\left| K_{n}\left( x-t\right) \right| d\mu \left( t\right) \leq
\frac{cM_{k}}{M_{N}},\qquad \text{when }n\geq M_{N}.
\end{equation*}

\begin{lemma}
Let $x\in I_{N}\left( x_{k}e_{k}+x_{l}e_{l}\right) ,$ $\ 1\leq x_{k}\leq
m_{k}-1,$ $1\leq x_{l}\leq m_{l}-1,$ $\ k=0,...,N-2,$ $l=k+1,...,N-1.$ Then
\end{lemma}

\begin{equation*}
\int_{I_{N}}\text{ }\underset{j=M_{N}+1}{\overset{n}{\sum }}\frac{\left|
K_{j}\left( x-t\right) \right| }{j+1}d\mu \left( t\right) \leq \frac{%
cM_{k}M_{l}}{M_{N}^{2}}.
\end{equation*}

Let $x\in I_{N}\left( x_{k}e_{k}\right) ,$ $1\leq x_{k}\leq m_{k}-1,$ $\
k=0,...,N-1.$ Then

\begin{equation*}
\int_{I_{N}}\underset{j=M_{N}+1}{\overset{n}{\sum }}\frac{\left| K_{j}\left(
x-t\right) \right| }{j+1}d\mu \left( t\right) \leq \frac{cM_{k}}{M_{N}}l_{n}.
\end{equation*}

\textbf{Proof. }Let $x\in I_{N}\left( x_{k}e_{k}+x_{l}e_{l}\right) ,$ $\
1\leq x_{k}\leq m_{k}-1,$ $1\leq x_{l}\leq m_{l}-1,$ $\ k=0,...,N-2,$ $%
l=k+1,...,N-1.$ Using Lemma 3 we have

\begin{eqnarray}
&&\int_{I_{N}}\underset{j=M_{N}+1}{\overset{n}{\sum }}\frac{\left\vert
K_{j}\left( x-t\right) \right\vert }{j+1}d\mu \left( t\right) \leq \underset{%
j=M_{N}+1}{\overset{n}{\sum }}\frac{cM_{k}M_{l}}{\left( j+1\right) jM_{N}}
\label{12a} \\
&\leq &\frac{cM_{k}M_{l}}{M_{N}}\underset{j=M_{N}+1}{\overset{\infty }{\sum }%
}\left( \frac{1}{j}-\frac{1}{j+1}\right) \leq \frac{cM_{k}M_{l}}{M_{N}^{2}}.
\notag
\end{eqnarray}

Let $x\in I_{N}\left( x_{k}e_{k}\right) ,$ $1\leq x_{k}\leq m_{k}-1,$ $\
k=0,...,N-1.$ Then

\begin{equation}
\int_{I_{N}}\underset{j=M_{N}+1}{\overset{n}{\sum }}\frac{\left\vert
K_{j}\left( x-t\right) \right\vert }{j+1}d\mu \left( t\right) \leq \underset{%
j=M_{N}+1}{\overset{n}{\sum }}\frac{cM_{k}}{\left( j+1\right) M_{N}}\leq
\frac{cM_{k}}{M_{N}}l_{n}.  \label{12b}
\end{equation}

Combining (\ref{12a}) and (\ref{12b}) we complete the proof of Lemma 4.

\section{Proof of the Theorems}

\textbf{Proof of theorem 1.} \textbf{a) }Using Abel transformation we obtain

\begin{equation}
R_{n}f=\frac{1}{l_{n}}\overset{n-1}{\underset{j=1}{\sum }}\frac{\sigma _{j}f%
}{j+1}+\frac{\sigma _{n}f}{l_{n}}.  \label{24}
\end{equation}%
Consequently,

\begin{equation}
R^{\ast }f\leq c\sigma ^{\ast }f.  \label{24a}
\end{equation}

Using Theorem W and (\ref{24a}) we conclude that $R^{\ast }$ is bounded from
the martingale Hardy space $H_{1/2}$ to the space $L_{1/2,\infty }$.

\textbf{Proof of theorem 2. }From (\ref{24}) for the kernels of Riesz`s
logarithmic means we have

\begin{equation}
L_{n}=\frac{1}{l_{n}}\overset{n-1}{\underset{j=1}{\sum }}\frac{K_{j}}{j+1}+%
\frac{K_{n}}{l_{n}}.  \label{25}
\end{equation}

By Lemma 1, the proof of theorem 2 will be complete, if we show that

\begin{equation*}
\int\limits_{\overset{-}{I}}\left\vert \overset{\sim }{R}^{\ast
}a\right\vert ^{1/2}d\mu \leq c<\infty ,
\end{equation*}%
for every 1/2-atom $a,$ where $I$ denotes the support of the atom.

Let $a$ be an arbitrary 1/2-atom with support$\ I$ and $\mu \left( I\right)
=M_{N}^{-1}.$ We may assume that $I=I_{N}.$ It is easy to see that $%
R_{n}\left( a\right) =\sigma _{n}\left( a\right) =0,$ when $n\leq M_{N}$.
Therefore we suppose that $n>M_{N}.$

Since $\left\Vert a\right\Vert _{\infty }\leq cM_{N}^{2}$ if we apply (\ref%
{25}) we can write
\begin{eqnarray}
&&\frac{\left\vert R_{n}a\left( x\right) \right\vert }{\log \left(
n+1\right) }=\frac{1}{\log \left( n+1\right) }\int_{I_{N}}\left\vert a\left(
t\right) \right\vert \left\vert L_{n}\left( x-t\right) \right\vert d\mu
\left( t\right)  \label{26} \\
&\leq &\frac{\left\Vert a\right\Vert _{\infty }}{\log \left( n+1\right) }%
\int_{I_{N}}\left\vert L_{n}\left( x-t\right) \right\vert d\mu \left(
t\right)  \notag \\
&\leq &\frac{cM_{N}^{2}}{\log \left( n+1\right) l_{n}}\underset{I_{N}}{\int }%
\text{ }\underset{j=M_{N}+1}{\overset{n-1}{\sum }}\frac{\left\vert
K_{j}\left( x-t\right) \right\vert }{j+1}d\mu \left( t\right)  \notag \\
&&+\frac{cM_{N}^{2}}{\log \left( n+1\right) l_{n}}\int_{I_{N}}\left\vert
K_{n}\left( x-t\right) \right\vert d\mu \left( t\right) .  \notag
\end{eqnarray}

Let $x\in I_{N}\left( x_{k}e_{k}+x_{l}e_{l}\right) ,$ $\ 1\leq x_{k}\leq
m_{k}-1,$ $1\leq x_{l}\leq m_{l}-1,$ $k=0,...,N-2,$ $l=k+1,...,N-1.$ From
Lemmas 3 and 4 we have

\begin{equation}
\frac{\left| R_{n}\left( a\right) \right| }{\log \left( n+1\right) }\leq
\frac{cM_{l}M_{k}}{N^{2}}.  \label{29}
\end{equation}

Let $x\in I_{N}\left( x_{k}e_{k}\right) ,$ $1\leq x_{k}\leq m_{k}-1,$ $%
k=0,...,N-1.$ Applying Lemmas 3 and 4 we have

\begin{equation}
\frac{\left\vert R_{n}a\left( x\right) \right\vert }{\log \left( n+1\right) }%
\leq \frac{M_{N}M_{k}}{N}\leq cM_{N}M_{k}.  \label{32}
\end{equation}

Combining (\ref{3}), (\ref{29}) and (\ref{32}) we get
\begin{eqnarray*}
&&\int_{\overline{I_{N}}}\left\vert \overset{\sim }{R}^{\ast }a\left(
x\right) \right\vert ^{1/2}d\mu \left( x\right) \\
&=&\overset{N-2}{\underset{k=0}{\sum }}\overset{m_{k}-1}{\underset{x_{k}=1}{%
\sum }}\overset{N-1}{\underset{l=k+1}{\sum }}\overset{m_{l}-1}{\underset{%
x_{l}=1}{\sum }}\int_{I_{l+1}\left( x_{k}e_{k}+x_{l}e_{l}\right) }\left\vert
\overset{\sim }{R}^{\ast }a\left( x\right) \right\vert ^{1/2}d\mu \left(
x\right) \\
&&+\overset{N-1}{\underset{k=0}{\sum }}\overset{m_{k}-1}{\underset{x_{k}=1}{%
\sum }}\int_{I_{N}\left( x_{k}e_{k}\right) }\left\vert \overset{\sim }{R}%
^{\ast }a\left( x\right) \right\vert ^{1/2}d\mu \left( x\right) \\
&\leq &c\overset{N-2}{\underset{k=0}{\sum }}\overset{N-1}{\underset{l=k+1}{%
\sum }}\frac{1}{M_{l}}\frac{\sqrt{M_{l}M_{k}}}{N}+c\overset{N-1}{\underset{%
k=0}{\sum }}\frac{1}{M_{N}}\sqrt{M_{N}M_{k}}\leq c<\infty .
\end{eqnarray*}

It completes the proof of first part of theorem 2.

\textbf{b)} Let$\ \left\{ \lambda _{k},\text{ }k\in \mathbf{P}_{+}\right\} $
be an increasing sequence of the positive integers, which saisfies condition
(\ref{55}). For every $\lambda _{k}$ there exists a positive integers$\
\left\{ n_{k},\text{ }k\in \mathbf{P}_{+}\right\} \subset \left\{ \lambda
_{k},\text{ }k\in \mathbf{P}_{+}\right\} ,$ such that
\begin{equation*}
\lim_{k\rightarrow \infty }\frac{n_{k}}{\varphi \left( M_{2n_{k}+1}\right) }%
=\infty .
\end{equation*}

Let

\begin{equation*}
f_{n_{k}}\left( x\right) =D_{M_{2n_{k}+1}}\left( x\right)
-D_{M_{_{2n_{k}}}}\left( x\right) .
\end{equation*}

It is evident
\begin{equation*}
\widehat{f}_{n_{k}}\left( i\right) =\left\{
\begin{array}{l}
\text{ }1,\text{ if }i=M_{_{2n_{k}}},...,M_{2n_{k}+1}-1, \\
\text{ }0,\text{otherwise}.%
\end{array}%
\right.
\end{equation*}%
We can write
\begin{equation}
S_{i}f_{n_{k}}\left( x\right) =\left\{
\begin{array}{l}
D_{i}\left( x\right) -D_{M_{_{2n_{k}}}}\left( x\right) ,\text{ \ if }%
i=M_{_{2n_{k}}},...,M_{2n_{k}+1}-1, \\
\text{ }f_{n_{k}}\left( x\right) ,\text{ \ if \ }i\geq M_{2n_{k}+1}, \\
0,\text{ \qquad otherwise.}%
\end{array}%
\right.  \label{33}
\end{equation}

From (\ref{4}) we get (see also \cite{tep2} and \cite{tep3})
\begin{equation}
\left\Vert f_{n_{k}}\left( x\right) \right\Vert _{H_{p}}=\left\Vert
f_{n_{k}}^{\ast }\left( x\right) \right\Vert _{p}\leq cM_{_{2n_{k}}}^{1-1/p}.
\label{34}
\end{equation}%
Let $q_{n_{k}}^{s}=M_{2n_{k}}+M_{2s},$ $s=0,...,n_{k}-1.$ By (\ref{33}) we
have%
\begin{eqnarray}
\frac{\left\vert R_{q_{n_{k}}^{s}}f_{n_{k}}\left( x\right) \right\vert }{%
\varphi \left( q_{n_{k}}^{s}\right) } &=&\frac{1}{\varphi \left(
q_{n_{k}}^{s}\right) l_{q_{n_{k}}^{s}}}\left\vert \overset{q_{n_{k}}^{s}}{%
\underset{j=M_{_{2n_{k}}}+1}{\sum }}\frac{S_{j}f_{n_{k}}\left( x\right) }{j}%
\right\vert  \label{35} \\
&=&\frac{1}{\varphi \left( q_{n_{k}}^{s}\right) l_{q_{n_{k}}^{s}}}\left\vert
\overset{q_{n_{k}}^{s}}{\underset{j=M_{_{2n_{k}}}+1}{\sum }}\frac{\left(
D_{j}\left( x\right) -D_{M_{_{2n_{k}}}}\left( x\right) \right) }{j}%
\right\vert  \notag \\
&=&\frac{1}{\varphi \left( q_{n_{k}}^{s}\right) l_{q_{n_{k}}^{s}}}\left\vert
\overset{M_{2s}}{\underset{j=1}{\sum }}\frac{\left(
D_{j+M_{_{2n_{k}}}}\left( x\right) -D_{M_{_{2n_{k}}}}\left( x\right) \right)
}{j+M_{_{2n_{k}}}}\right\vert .  \notag
\end{eqnarray}

Since

\begin{equation}
D_{j+M_{_{2n_{k}}}}\left( x\right) -D_{M_{_{2n_{k}}}}\left( x\right) =\psi
_{M_{_{2n_{k}}}\text{ }}D_{j}\left( x\right) ,\text{ \ \ }%
j=1,2,..,M_{_{2n_{k}}}-1.  \label{36}
\end{equation}%
we obtain

\begin{equation}
\frac{\left\vert R_{q_{n_{k}}^{s}}f_{n_{k}}\left( x\right) \right\vert }{%
\varphi \left( q_{n_{k}}^{s}\right) }=\frac{1}{\varphi \left(
q_{n_{k}}^{s}\right) l_{q_{n_{k}}^{s}}}\overset{M_{2s}}{\underset{j=1}{\sum }%
}\frac{\left\vert Dj\left( x\right) \right\vert }{j+M_{_{2n_{k}}}}.
\label{37}
\end{equation}

Let $x\in $ $I_{_{2s}}\backslash I_{_{2s+1}}.$ Then

\begin{eqnarray}
&&\frac{\left\vert R_{q_{n_{k}}^{s}}f_{n_{k}}\left( x\right) \right\vert }{%
\varphi \left( q_{n_{k}}^{s}\right) }\geq \frac{1}{\varphi \left(
q_{n_{k}}^{s}\right) l_{q_{n_{k}}^{s}}}\overset{M_{2s}}{\underset{j=0}{\sum }%
}\frac{j}{j+M_{_{2n_{k}}}}  \label{38} \\
&\geq &\frac{1}{\varphi \left( q_{n_{k}}^{s}\right) l_{q_{n_{k}}^{s}}}\frac{%
\overset{M_{2s}}{\underset{j=0}{\sum }}j}{2M_{_{2n_{k}}}}\geq \frac{%
cM_{2s}^{2}}{\varphi \left( q_{n_{k}}^{s}\right)
l_{q_{n_{k}}^{s}}M_{_{2n_{k}}}}.  \notag
\end{eqnarray}

Using (\ref{38}) we have
\begin{eqnarray*}
&&\int_{G_{m}}\left\vert \overset{\sim }{R}^{\ast }f\left( x\right)
\right\vert ^{1/2}d\mu \left( x\right) \\
&\geq &\text{ }\overset{n_{k}-1}{\underset{s=1}{\sum }}\int_{^{I_{_{2s}}%
\backslash I_{_{2s+1}}}}\left\vert \frac{R_{q_{n_{k}}^{s}}f\left( x\right) }{%
\varphi \left( q_{n_{k}}^{s}\right) }\right\vert ^{1/2}d\mu \left( x\right)
\geq c\overset{n_{k}-1}{\underset{s=1}{\sum }}\frac{M_{2s}}{\sqrt{\varphi
\left( q_{n_{k}}^{s}\right) l_{q_{n_{k}}^{s}}M_{_{2n_{k}}}}}\frac{1}{M_{2s}}
\\
&\geq &c\overset{n_{k}-1}{\underset{s=1}{\sum }}\frac{1}{\sqrt{\varphi
\left( M_{_{2n_{k}+1}}\right) l_{M_{_{2n_{k}+1}}}M_{_{2n_{k}}}}}\geq \frac{%
cn_{k}}{\sqrt{\varphi \left( M_{_{2n_{k}+1}}\right)
l_{M_{_{2n_{k}+1}}}M_{_{2n_{k}}}}}.
\end{eqnarray*}%
From (\ref{34}) we have
\begin{equation}
\frac{\left( \int_{G_{m}}\left\vert \overset{\sim }{R}^{\ast }f\left(
x\right) \right\vert ^{1/2}d\mu \left( x\right) \right) ^{2}}{\left\Vert
f_{n_{k}}\left( x\right) \right\Vert _{H_{1/2}}}\geq \frac{cn_{k}}{\varphi
\left( M_{_{2n_{k}+1}}\right) }\rightarrow \infty ,\text{ when }k\rightarrow
\infty .  \label{40}
\end{equation}

Theorem 2 is proved.

\textbf{Proof of theorem 3.} Let $0<p<1/2.$\textbf{\ }By Lemma 1, the proof
of theorem 3 will be complete, if we show that

\begin{equation*}
\int\limits_{\overset{-}{I}}\left\vert \overset{\sim }{R}_{p}^{\ast
}a\right\vert ^{p}d\mu \leq c_{p}<\infty ,\text{ }
\end{equation*}%
for every p-atom $a,$ where $I$ denotes the support of the atom$.$

Let $a$ be an arbitrary p-atom with support$\ I$ and $\mu \left( I\right)
=M_{N}^{-1}.$ We may assume that $I=I_{N}.$ It is easy to see that $%
R_{n}\left( a\right) =0,$ when $n\leq M_{N}$. Therefore we suppose that $%
n>M_{N}.$

Since $\left\| a\right\| _{\infty }\leq cM_{N}^{1/p}$ using (\ref{25}) we
can write

\begin{eqnarray}
&&\frac{\log \left( n+1\right) }{\left( n+1\right) ^{1/p-2}}\left\vert
R_{n}a\left( x\right) \right\vert  \label{41} \\
&\leq &\frac{\log \left( n+1\right) M_{N}^{1/p}}{\left( n+1\right)
^{1/p-2}l_{n}}\underset{I_{N}}{\int }\text{ }\underset{j=M_{N}+1}{\overset{%
n-1}{\sum }}\frac{\left\vert K_{j}\left( x-t\right) \right\vert }{j+1}d\mu
\left( t\right)  \notag \\
&&+\frac{\log \left( n+1\right) M_{N}^{1/p}}{\left( n+1\right) ^{1/p-2}l_{n}}%
\int_{I_{N}}\left\vert K_{n}\left( x-t\right) \right\vert d\mu \left(
t\right) .  \notag
\end{eqnarray}

Let $x\in I_{N}\left( x_{k}e_{k}+x_{l}e_{l}\right) ,$ $\ 1\leq x_{k}\leq
m_{k}-1,$ $1\leq x_{l}\leq m_{l}-1,$ $k=0,...,N-2,$ $l=k+1,...,N-1.$ From
Lemmas 3 and 4 when $n>M_{N}$ we obtain
\begin{equation}
\frac{\log \left( n+1\right) }{\left( n+1\right) ^{1/p-2}}\left\vert
R_{n}a\left( x\right) \right\vert \leq \text{ }c_{p}M_{l}M_{k}.  \label{42a}
\end{equation}

Let $x\in I_{N}\left( x_{k}e_{k}\right) ,$ $1\leq x_{k}\leq m_{k}-1,\ $\ $%
k=0,...,N-1.$ Applying Lemmas 3 and 4 we have

\begin{equation}
\frac{\log \left( n+1\right) }{\left( n+1\right) ^{1/p-2}}\left\vert
R_{n}a\left( x\right) \right\vert \leq cNM_{N}M_{k}.  \label{44}
\end{equation}

Combining (\ref{3}), (\ref{42a}) and (\ref{44}) we get

\begin{eqnarray*}
&&\int_{\overline{I_{N}}}\left\vert \overset{\sim }{R}_{p}^{\ast }a\left(
x\right) \right\vert ^{p}d\mu \left( x\right) \\
&=&\overset{N-2}{\underset{k=0}{\sum }}\overset{m_{k}-1}{\underset{x_{k}=1}{%
\sum }}\overset{N-1}{\underset{l=k+1}{\sum }}\overset{m_{l}-1}{\underset{%
x_{l}=1}{\sum }}\int_{I_{N}\left( x_{k}e_{k}+x_{l}e_{l}\right) }\left\vert
\overset{\sim }{R}_{p}^{\ast }a\left( x\right) \right\vert ^{p}d\mu \left(
x\right) \\
&&+\overset{N-1}{\underset{k=0}{\sum }}\overset{m_{k}-1}{\underset{x_{k}=1}{%
\sum }}\int_{I_{N}\left( x_{k}e_{k}\right) }\left\vert \overset{\sim }{R}%
_{p}^{\ast }a\left( x\right) \right\vert ^{p}d\mu \left( x\right) \\
&\leq &c_{p}\overset{N-2}{\underset{k=0}{\sum }}\overset{N-1}{\underset{l=k+1%
}{\sum }}\frac{1}{M_{l}}\left( M_{l}M_{k}\right) ^{p}+c_{p}\overset{N-1}{%
\underset{k=0}{\sum }}\frac{1}{M_{N}}\left( NM_{N}M_{k}\right) ^{p}\leq
c_{p}<\infty .
\end{eqnarray*}

Which complete the proof of first part of Theorem 2.

Let $0<p<1/2$ and$\ \left\{ \lambda _{k},\text{ }k\in \mathbf{P}_{+}\right\}
$ be an increasing sequence of the positive integers, which satisfies
condition (\ref{66}). It is evident that for every $\lambda _{k}$ there
exists a positive integers $\ \left\{ n_{k},\text{ }k\in \mathbf{P}%
_{+}\right\} \subset \left\{ \lambda _{k},\text{ }k\in \mathbf{P}%
_{+}\right\} ,$ such that
\begin{equation*}
\lim_{k\rightarrow \infty }\frac{\left( M_{2n_{k}}+1\right) ^{1/p-2}}{%
\varphi \left( M_{2n_{k}}+1\right) \log \left( M_{2n_{k}}+1\right) }=\infty .
\end{equation*}

Combining (\ref{35}-\ref{38}) we have

\begin{equation*}
\frac{\left\vert R_{M_{2n_{k}}+1}f_{n_{k}}\left( x\right) \right\vert }{%
\varphi \left( M_{2n_{k}}+1\right) }=\frac{\left\vert
R_{q_{n_{k}}^{0}}f\left( x\right) \right\vert }{\varphi \left(
q_{n_{k}}^{0}\right) }\geq \frac{c}{\varphi \left( M_{_{2n_{k}}}+1\right)
l_{M_{_{2n_{k}}}+1}\left( M_{_{2n_{k}}}+1\right) },
\end{equation*}%
for $x\in I_{0}\backslash I_{1}=G_{m}\backslash I_{1}.$

From (\ref{34}) we get

\begin{eqnarray*}
&&\frac{\frac{c}{\varphi \left( M_{_{2n_{k}}}+1\right)
l_{M_{_{2n_{k}}}+1}\left( M_{_{2n_{k}}}+1\right) }\mu \left\{ x\in
G_{m}:\left\vert \overset{\sim }{R}_{p}^{\ast }f_{n_{k}}\left( x\right)
\right\vert \geq \frac{c}{\varphi \left( M_{_{2n_{k}}}+1\right)
l_{M_{_{2n_{k}}}+1}\left( M_{_{2n_{k}}}+1\right) }\right\} ^{1/p}}{%
\left\Vert f_{n_{k}}\left( x\right) \right\Vert _{H_{p}}} \\
&\geq &\frac{c\left( M_{2n_{k}}+1\right) ^{1/p-2}}{\varphi \left(
M_{2n_{k}}+1\right) \log \left( M_{2n_{k}}+1\right) }\rightarrow \infty ,%
\text{ when }k\rightarrow \infty .
\end{eqnarray*}

Which complete the proof of theorem 3.

\textbf{Acknowledgment: }The author would like to thank the referee for
helpful suggestions.

\end{document}